\documentclass[titlepage,fleqn,11pt]{article}

\usepackage[margin=1.2in]{geometry}
\usepackage{latexsym}
\usepackage{amsfonts}
\usepackage[leqno]{amsmath}

\makeatletter
\newcommand{\leqnomode}{\tagsleft@true}
\newcommand{\reqnomode}{\tagsleft@false}
\makeatother

\usepackage{tikz}
\usepackage{float}
\usepackage{url}
\usetikzlibrary{arrows}
\usetikzlibrary{decorations.markings}
\usepackage[symbol]{footmisc}

\usepackage{xspace}
\usepackage[color=Olive]{todonotes}

\usepackage[colorlinks=true,allcolors=blue]{hyperref}

\def\longbox#1{\parbox{0.85\textwidth}{#1}}

\newcommand{\lem}[1]{\vspace*{-0.4cm}
  \begin{equation}\vspace*{-0.4cm}
    \longbox{\emph{#1}}
\end{equation}}

\newcommand\blackslug{\hbox{\hskip 1pt \vrule width 4pt height 8pt depth 1.5pt
        \hskip 1pt}}
\newcommand\bbox{\hfill \quad \blackslug \bigbreak}

\newcommand*\samethanks[1][\value{footnote}]{\footnotemark[#1]}

\usepackage{authblk}

\title{On Aharoni's rainbow generalization of the Caccetta-H\"{a}ggkvist conjecture}

\author[1]{Patrick Hompe}
\author[2]{Petra Pelik\'{a}nov\'{a}\thanks{This paper is partially based on research performed at the DIMACS REU 2019, which has been supported by the H2020-MSCA-RISE project CoSP- GA No. 823748.}}
\author[2]{Aneta Pokorn\'{a}\samethanks[1]\thanks{Supported by GAUK 1277018.}}
\author[1]{Sophie Spirkl\thanks{Part of this research was performed while the third author was at Rutgers University.

This material is based upon work supported by the National Science Foundation under Award No. DMS-1802201.

We acknowledge the support of the Natural Sciences and Engineering Research Council of Canada (NSERC), [funding reference number RGPIN-2020-03912]. Cette recherche a été financée par le Conseil de recherches en sciences naturelles et en génie du Canada (CRSNG), [numéro de référence RGPIN-2020-03912].
}}

\affil[1]{University of Waterloo}
\affil[2]{Charles University}

\date{\today}

\newtheorem{thm}{Theorem}[section]
\newtheorem{cor}[thm]{Corollary}
\newtheorem{conjecture}[thm]{Conjecture}

\newcommand{\Proof}{\noindent{\bf Proof.}\ \ }

\begin{document}
\maketitle
\begin{abstract}
For a digraph $G$ and $v \in V(G)$, let $\delta^+(v)$ be the number of out-neighbors of $v$ in $G$. The Caccetta-H\"{a}ggkvist conjecture states that for all $k \ge 1$, if $G$ is a digraph with $n = |V(G)|$ such that $\delta^+(v) \ge k$ for all $v \in V(G)$, then G contains a directed cycle of length at most $\lceil n/k \rceil$. In [\ref{aharoni}], Aharoni proposes a generalization of this conjecture, that a simple edge-colored graph on $n$ vertices with $n$ color classes, each of size $k$, has a rainbow cycle of length at most $\lceil n/k \rceil$. In this paper, we prove this conjecture if each color class has size ${\Omega}(k \log k)$.
\end{abstract}

\section{Introduction and preliminaries}
A graph or digraph is \emph{simple} if there are no loops or parallel edges. For a simple digraph $G$ and a vertex
$v \in V(G)$, let $\delta^+(v)$ denote the number of out-neighbors of $v$ in $G$. A famous conjecture in graph theory is the following, due to Caccetta and H\"{a}ggkvist [\ref{ch}]:

\begin{conjecture}[Caccetta-H\"{a}ggkvist]\label{ch_thm}
Suppose $n,k$ are positive integers, and let $G$ be a simple digraph on $n$ vertices with $\delta^+(v) \ge k$ for all $v \in V(G)$; then $G$ contains a directed cycle of length at most $\lceil n/k \rceil$.
\end{conjecture}

For a graph $G$ and a function $c: E(G) \rightarrow \{1, \dots, |V(G)|\}$, a \emph{rainbow cycle (with respect to $c$)} is a cycle $C$ in $G$ such that for all $e, f \in E(C)$ with $e \neq f$, we have $c(e) \neq c(f)$. We will refer to $c$ as a \emph{coloring} of the edges of $G$.\footnote{Note that $c$ is not required to be a proper edge-coloring.} We say that $c$ has \emph{color classes of size at least $k$} for $k \in \mathbb{N}$ if $|c^{-1}(i)| \geq k$ for all $i \in \{1, \dots, |V(G)|\}$. 

In [\ref{aharoni}], Aharoni proposes a generalization of Conjecture \ref{ch_thm}:

\begin{conjecture}[\cite{aharoni}]\label{aharoni_thm}
Let $n,k$ be positive integers, and let $G$ be a simple graph on $n$ vertices. Let $c : E(G) \rightarrow \{1, \dots, n\}$ be a coloring of the edges of $G$ with color classes of size at least $k$; then $G$ has a rainbow cycle of length at most $\lceil n/k \rceil$.
\end{conjecture}

In a recent paper, Devos et al.\ [\ref{devos}] prove that Conjecture \ref{aharoni_thm} is true for $k=2$:

\begin{thm}[\cite{devos}]\label{devos_thm}
Let $G$ be a simple graph on $n$ vertices, and let $c$ be a coloring of the edges of $G$ with color classes of size at least $2$; then there exists a rainbow cycle of length at most $\lceil n/2 \rceil$.
\end{thm}

We also make use of the following results due to Bollob\'{a}s and Szemer\'{e}di [\ref{bollobas}] and Shen [\ref{shen}], respectively. The first deals with the girth of a simple graph, while the second is an approximate result for Conjecture \ref{ch_thm}. In this paper, $\log$ denotes the logarithm with base $2$.

\begin{thm}[\cite{bollobas}]\label{bollobas_thm}
For all $n \ge 4$ and $k \ge 2$, if $G$ is a simple graph on $n$ vertices with $n+k$ edges, then $G$ contains a cycle of length at most $$\frac{2(n+k)}{3k}(\log k+\log \log k + 4).$$
\end{thm}

\begin{thm}[\cite{shen}]\label{shen_thm}
Let $G$ be a simple digraph with $\delta^+(v) \ge k$ for all $v \in V(G)$. Then $G$ contains a directed cycle of length at most $\lceil n/k \rceil + 73$.
\end{thm}

\section{Main result}
Our main result is the following: 
\begin{thm}\label{main}
Let $k > 1$ be an integer, and let $G$ be a graph. Let $c$ be a coloring of the edges of $G$ with color classes of size at least $301k \log k$. Then $G$ contains a rainbow cycle of length at most $\lceil n/k \rceil$.
\end{thm}
\Proof
We proceed by induction on the number of vertices. Let $f(k) = 7k \log k$, and let $G$ be a graph on $n$ vertices. Let $c$ be a coloring of the edges of $G$ with color classes of size at least $43f(k)$. Suppose for a contradiction that there is no rainbow cycle of length at most $\lceil n/k \rceil$. Note that $G$ has at least $43f(k)n$ edges, and therefore, $n > 43f(k)$.

For $v \in V(G)$ and $i \in \{1, \dots, n\}$, we say that $i$ \emph{is dominant at} $v$ if $v$ is incident with at least $7f(k)$ edges $e$ such that $c(e) = i$. We call a vertex $v \in V(G)$ \textit{color-dominating} if there exists $i \in \{1, \dots, n\}$ such that $i$ is dominant at $v$. We call a color $i \in \{1, \dots, n\}$ \textit{vertex-dominating} if there exists a vertex $v \in V(G)$ such that $i$ is dominant at $v$. Let us say that $H \subseteq V(G)$ is \textit{nice} if
\begin{itemize}
    \item for every vertex-dominating color $i \in \{1, \dots, n\}$, there is a vertex $v \in V(G) \setminus H$ such that $i$ is dominant at $v$; and
    \item there are at most $|H|$ colors $i \in \{1, \dots, n\}$ such that $i$ is not vertex-dominating and for all $e \in c^{-1}(i)$, at least one end of $e$ is in $H$.
\end{itemize} 
 
\lem{\label{lemma1}
If there is a nice set $H \subseteq V(G)$ with $6f(k) \le |H| < n$, then there is a nice set $H' \subseteq V(G)$ with $|H'| = \lceil 6f(k) \rceil$.}

We remove vertices from $H$ one-by-one such that the remaining set is nice. Suppose that we have removed $j \ge 0$ vertices from $H$, leaving a nice set $H_j$ with $|H_j| > \lceil 6f(k) \rceil$. Let $C_j$ be the set of colors $i \in \{1, \dots, n\}$ which are not vertex-dominating and also do not have an edge $e$ with $c(e) = i$ such that both ends of $e$ are in $V(G) \setminus H_j$. From the definition of a nice set, we know $|C_j| \le |H_j|$. If $|C_j| < |H_j|$, then removing any vertex from $H_j$ gives a smaller nice set. So, we may assume that $|C_j| = |H_j|$. If there is a color $i$ in $C_j$ and an edge $e = uv \in c^{-1}(i)$ with $v \in H_j$ and $u \in G \setminus H_j$, then $H_j \setminus \{v\}$ is nice. If there is no such $i \in C_j$, then for every color $i \in C_j$, all edges in $c^{-1}(i)$ have both their ends in $H_j$. Now applying induction to the subgraph of $G$ with vertex set $H_j$ and edge set $c^{-1}(C_j)$ gives a rainbow cycle of length at most $\lceil n/k \rceil$ in $G$, a contradiction. This proves \ref{lemma1}. 

\lem{\label{lemma2}
There is a nice set $H' \subseteq V(G)$ with $|H'| = \lceil 6f(k) \rceil$.}

For each vertex-dominating color $i$, we pick a vertex $v_i$ such that $i$ is dominant at $v_i$, and let $S$ be the set of these vertices $v_i$. Let $H=V(G) \setminus S$. Note that $H$ is nice;
thus by \ref{lemma1}, we may assume that either $|H| < 6f(k)$ or $|H| = n$. 

We first consider the case when $|H| = n$. Since $43f(k) \ge 2$, Theorem \ref{devos_thm} guarantees the existence of a rainbow cycle $K$ of length at most $n/2+1$ in $G$. Let $H' = V(G) \setminus  V(K)$. Then $H'$ is nice, and $n > |H'| \ge n/2 - 1 \ge 6f(k)$; so by  \ref{lemma1}, $G$ contains a nice set of size $\lceil 6f(k) \rceil$.

Now we may assume that $|H| < 6f(k)$. We construct a digraph $G'$ with $V(G') = S$, and for all $i, j$ with $v_i, v_j \in S$, there is an arc $v_i \rightarrow v_j$ if $v_iv_j \in E(G)$ and $c(v_iv_j) = i$. Every vertex $v_i$ is incident with at least $7f(k)$ edges $e$ with $c(e) = i$, and since $|H| < 6f(k)$, there are at least $f(k)$ edges $e = v_iu$ with $c(e) = i$ and $u \in S$. Therefore, $\delta^+(G') \geq f(k)$. 

Now, we claim $n/f(k) + 74 \le n/k$, which is equivalent to $74kf(k) \le n(f(k)-k)$. Since $k \ge 2$, we have $\log(k) \ge 117/301$, and thus $74kf(k) \le 43f(k)(f(k)-k) \le n(f(k)-k)$, as claimed.

Then, by applying Theorem \ref{shen_thm} to $G'$ we obtain a directed cycle $K$ of length at most $\lceil n/f(k) \rceil + 73 \le \lceil n/k \rceil$ in $G'$. The edges of $G$ that correspond to arcs of $K$ form a rainbow cycle of length at most $\lceil n/k \rceil$ in $G$, a contradiction. This proves \ref{lemma2}.\bbox{}
\lem{\label{lemma3}
Let $H \subseteq V(G)$ be a nice set with $|H| = \lceil 6f(k) \rceil$. Then there exists $H' \subseteq H$ such that $|H'| \geq  \lceil 2f(k) \rceil$ and such that for at least $n - \lceil{f(k)}\rceil + 1$ colors $i$, at least one edge $e \in c^{-1}(i)$ has both ends in $V(G) \setminus H'$.}

Let $C$ be the set of colors $i$ which are not vertex-dominating and for which no edge of $c^{-1}(i)$ has both ends in $V(G) \setminus H$. Since $H$ is nice, it follows that $|C| \le |H| = \lceil 6f(k) \rceil$. Let $D \subseteq C$ be the set of colors $i \in C$ such that there is a vertex $v \in H$ which is incident with all edges in $c^{-1}(i)$ that have one end in $H$ and the other in $V(G) \setminus H$. We claim that $|D| \le \lceil f(k) \rceil - 1$. Indeed, for each color $i \in D$, there are at least $\lceil{36f(k)}\rceil$ edges in $c^{-1}(i)$ with both ends in $H$ since $i$ is not vertex-dominating. If $|D| > \lceil f(k)\rceil - 1$, then we obtain more than $(f(k)-1)(36f(k))$ edges with both ends in $H$. Now, since $k \ge 2$, we have $f(k) \ge 72/23$, and it follows that:
$$(f(k)-1)(36f(k)) \ge \frac{49f(k)^2}{2} \ge \frac{(6f(k)+1)^2}{2} \ge \frac{|H|^2}{2}$$
which gives a contradiction. Thus, $|D| \le \lceil f(k)\rceil - 1$.

Next, we claim there exists $H' \subseteq H$ such that $|H'| = \lceil 2f(k) \rceil$ and such that for all $i \in \{1, \dots, n\} \setminus D$, there is an edge $e \in c^{-1}(i)$ with both ends in $V(G) \setminus H'$. To see this, we construct a graph $J$ with vertex set $H$ and the following set of edges. For each $i \in C \setminus D$, we choose two vertices $v_1^i, v_2^i \in H$, each incident with an edge in $c^{-1}(i)$ whose other end is in $V(G) \setminus H$; we know from the definition of $D$ that this is possible. Now,  the graph $J$ has $|H|$ vertices and at most $|H|$ edges, and so $J$ has a stable set $H' \subseteq V(J)$ of size at least $|V(J)|/3 \geq 2f(k)$; and so $|H'| \geq \lceil{2f(k)}\rceil$. 

Now, for every color $i \in C \setminus D$, $V(G) \setminus H'$ contains at least one of $v_1^i, v_2^i$, and therefore, there is an edge in $c^{-1}(i)$ with both ends in $V(G) \setminus H'$.  Moreover, for every $i \in \{1, \dots, n\} \setminus C$, either $i$ dominates a vertex $v$ in $V(G) \setminus H \subseteq V(G) \setminus H'$ (and so, since $|H'| < 7f(k)$, there is an edge in $c^{-1}(i)$ incident with $v$ whose other end is not in $H'$); or there is an edge in $c^{-1}(i)$ with both ends in $V(G) \setminus H \subseteq V(G) \setminus H'$.  Thus, for at least $n - |D| \geq n - \lceil f(k) \rceil + 1$ colors $i$, at least one edge in $c^{-1}(i)$ has both ends in $V(G) \setminus H'$. This proves \ref{lemma3}.

\medskip

By combining \ref{lemma2} and \ref{lemma3}, we conclude that there exists $H' \subseteq V(G)$ with $|H'| \geq \lceil{2f(k)}\rceil$, and such that for at least $n - \lceil f(k) \rceil + 1$ colors $i$, at least one edge in $c^{-1}(i)$ has both ends in $V(G) \setminus H'$. Let $H''$ be a subgraph of $G$ with vertex set $V(G) \setminus H'$, obtained by taking exactly one edge in $c^{-1}(i)$ with both ends in $V(G) \setminus H'$ for all $i \in \{1, \dots, n\}$ which have such an edge. It follows that $|E(H'')| \geq |V(H'')| + \lceil f(k) \rceil$.

Now, we claim that $\frac{2(n+f(k))}{3(f(k))}(\log \log (f(k)) + \log (f(k)) + 4) \leq \frac{n}{k}$. Using $f(k) < n/43$, it suffices to show:
$$\frac{88(\log \log (f(k)) + \log (f(k)) + 4)}{129} \le 7\log(k)$$

Let $g(k) = 7\log(k) - \frac{88}{129}(\log \log (f(k)) + \log (f(k)) + 4)$. We have that $g(2) > 0$, and for $k \ge 2$ we have:
$$f(k)g'(k) \ln(2) = 49\log(k) - \frac{88}{129}f'(k)\left(\frac{1}{\log(f(k))\ln(2)}+1\right) > 0$$
since for $k \ge 2$ we have:
$$f'(k)\left(\frac{1}{\log(f(k))\ln(2)}+1\right) < (7+7\log(k))(3) \le 49\log(k)$$
So $g'(k) > 0$ for $k \ge 2$, and it follows that $g(k) \ge 0$ for $k \ge 2$, as desired.

Then, Theorem \ref{bollobas_thm} gives a rainbow cycle of length at most  $\frac{2(n+f(k))}{3(f(k))}(\log \log (f(k)) + \log (f(k)) + 4) \leq \lceil \frac{n}{k}\rceil$, a contradiction.
This proves Theorem \ref{main}.\bbox{}

We have an immediate corollary which gives us a result for the case of $\Omega(n \log n)$ color classes each of size $k$:

\begin{cor}\label{maincor}
Let $k$ be a positive integer and let $G$ be a simple graph on $n$ vertices. Let $c : E(G) \rightarrow \{1, \dots, t\}$  with $t \geq 303n \log n$, and with $|c^{-1}(i)| \geq k$ for all $i \in \{1, \dots, t\}$. Then $G$ contains a rainbow cycle in $G$ of length at most $\lceil n/k \rceil$.
\end{cor}
\Proof
Note that $t \geq 303n \log n \ge 303n \log k$. Since $303n\log k \ge n \lceil 301 \log k \rceil$, we can partition $\{1, \dots, t\}$ into $n$ parts, each of size at least $\lceil 301 \log k \rceil$; that is, there is a function $f : \{1, \dots, t\} \rightarrow \{1, \dots, n\}$ such that $|f^{-1}(i)| \geq \lceil 301 \log k \rceil$ for all $i \in \{1, \dots, n\}$.  By Theorem \ref{main}, applied to $G$ and $f \circ c$, we obtain a rainbow cycle of length at most $\lceil n/k \rceil$ in $G$ with respect to $f \circ c$, which is also rainbow with respect to $c$. This proves Corollary \ref{maincor}.\bbox{}
\section*{Acknowledgments}
This paper is partially based on research performed at the DIMACS REU 2019, which has been supported by the H2020-MSCA-RISE project CoSP- GA No.  823748. The third author was supported by GAUK 1277018. Part of this research was performed while the third author was at Rutgers University. This material is based upon work supported by the National Science Foundation under Award No.DMS-1802201. We acknowledge the support of the Natural Sciences and Engineering Research Council of Canada (NSERC), [funding reference number RGPIN-2020-03912]. Cette recherche a été financée par le Conseil de recherches en sciences naturelles et en génie du Canada (CRSNG), [numéro de référence RGPIN-2020-03912].

\end{document}